\documentclass{elsart}

\usepackage{graphicx}
\usepackage{subfigure}
\usepackage{subeqnarray}

\newcommand {\Inv}[1]{ \frac{1}{ #1 } }
\newcommand {\Dif}[2]{ \frac{\partial #1}{\partial #2} }
\newcommand {\Difd}[2]{ \frac{\partial^2 #1}{\partial #2 ^2} }

\newcommand {\Laplace}{ \nabla^2}
\newcommand {\V}[1]{{\bf #1}}

\newcommand {\veck}{ \V{k} }
\newcommand {\vecx}{ \V{x} }

\newcommand {\F}[1]{ \tilde{#1} }

\begin{document}

\begin{frontmatter}

\title{A Stable Explicit Scheme for Solving Inhomogeneous Constant
  Coefficients Differential Equation using Green's Function}
\thanks[email]{E-mail address: hiroshi.abe@na-net.ornl.gov}
\author{Hiroshi Abe\thanksref{email}}
%\address{3-56-10 Kamiarai, Tokorozawa-shi, Saitama-ken, Japan}

\begin{abstract}
  A numerical explicit method to evaluates transient solutions of
  linear partial differential inhomogeneous equation with constant
  coefficients is proposed. A general form of the scheme for a
  specific linear inhomogeneous equation is shown. The method is
  applied to the wave equation and the diffuse equation and
  is investigated by simulating simple models. The numerical solutions of the
  proposed method show good agreement to the exact solutions. Comparing
  with explicit FDM, FDM shows the instability by the violation of CFL
  condition whereas the proposed method is always stable irrespective
  of any time step width.
\end{abstract}


\end{frontmatter}

\section{Introduction}
\label{sec:intro}
Many numerical method have been proposed to solve partial differential
equations.
Many of the proposed method is for spatial modeling.
In spite of the variety of the spatial modeling, temporal modelings may
not have been proposed so many.

As for the temporal modeling, 
when you solve a harmonic oscillating model in time, the time
differentiation can be replace with $i\omega$.
When you solve a transient model in time, the time differentiation
is replaced with finite differentiation.
The finite differentiation can be said to be the most popular method for
the transient numerical solution.

The finite differentiation can be classified roughly into two
categories.
One is the explicit scheme and the rest is the implicit
scheme \cite{anderson}. 
The explicit scheme is easy to program and the amounts of computation
for advancing one time step is smaller than the implicit method.
One of the bad things is the time step width must be small enough to
satisfy a condition otherwise numerical instability may occur.
As for the implicit scheme, the computation is generally large for one
time step but it is always stable for any time step width.
From the user's point of view the explicit method may be comfortable for
programming.

In the meantime, Green's functions are well used in theoretical physics
but are not used in the numerical transient problems.
This is because if you adopt Green's function to solve a certain
transient problem you have to integrate all over the history in time and
space for each time step.
Although the Green's function may not be suitable for numerical simulation,
the Green's function is the response function of the basic equation
against Dirac's $\delta$ function or another functions.
This means that the Green's function gives the exact time variation of
the system against a certain stimulation.

The reason of the instability of the explicit Finite Difference Method
(FDM) is caused from the 
difficulty of the approximation of $\partial/\partial t$ with the
several order of polynomials of $\Delta t$.
There is no physical or mathematical inevitability for the approximation.
Whereas the Green's function is the exact solution of the basic equation
so the time variation predicted by the Green's function has the
inevitability.
It is a big motivation that the benefit of the Green's function is
to be integrated into the explicit time advancing to formulate a stable
explicit scheme irrespective of any time step width.

The present paper describes a new scheme for transient solutions using
Green's function to solve inhomogeneous linear differential equation
with constant coefficients.
Although the scheme is a explicit scheme, it is always stable
irrespective of any time step width different from the explicit FDM.

The preliminary works was applied to electromagnetics and diffuse
equations \cite{habe:ipsj}\cite{habe:cefc93}\cite{habe:jsiam93}.
In the present paper the general formulation of the scheme is
described. 

In the paper, generalized basic formulation is derived in section
\ref{sec:basic}. 
Section \ref{sec:tgm} describes specific formulations for wave and
diffuse equations.
Numerical simulation using the proposed method is investigated in
section \ref{sec:numerical} to show the effectiveness of the presented
method.
The numerical accuracy and stability is discussed in section
\ref{sec:discuss}. 

\section{Basic Formulation}
\label{sec:basic}

Suppose a partial differential equation,
\begin{equation}
  {\cal D}_t u + {\cal D}_x u = S(x,t),
  \label{eq:the_eq}
\end{equation}
where ${\cal D}_t, {\cal D}_x$ and $u$ are partial differential operator with
respect to time $t$, space $\vecx = (x,y,z)$ and a certain physical value,
respectively.
The partial differential operators ${\cal D}_t, {\cal D}_x$ are represented as,
\begin{eqnarray}
  {\cal D}_t & = & \sum_n a_n \frac{\partial^n}{\partial t^n}, \label{eq:t_dif} \\
  {\cal D}_x & = & \sum_{l,m,n} b_{lmn} \frac{\partial^l}{\partial x^l}
  \frac{\partial^m}{\partial y^m} \frac{\partial^n}{\partial z^n}, \label{eq:x_dif}
\end{eqnarray}
where $a_n, b_{lmn}$ are constant values.
The authors are now interested in solving Eq.(\ref{eq:the_eq}) to obtain the
transient solution $u$ for the given source function $S(\vecx,t)$.

To solve the Eq.(\ref{eq:the_eq}) by the proposed method, Fourier
transformation is applied to the basic equation at first.
A physical value can be expressed with integral form of Fourier
components as,
\begin{equation}
  A(\vecx,t) = \Inv{(2\pi)^{3/2}} \int_{-\infty}^{\infty} \F{A}(\veck,t)
  e^{i\veck\vecx} d\veck,
  \label{eq:fourier}
\end{equation}
where $\veck = (k_x,k_y,k_z)$ is the wave number.
Substituting Eq.(\ref{eq:fourier}) into Eq.(\ref{eq:the_eq}) we have,
\begin{equation}
  {\cal D}_t \F{u}(\veck,t) + {\cal K} \F{u}(\veck,t) = \F{S}(k,t),
  \label{eq:basic_eq}
\end{equation}
where
\begin{equation}
  {\cal K} = \sum_{l,m,n} b_{lmn} (ik_x)^l(ik_y)^m(ik_z)^n,
\end{equation}
$\F{u}$ and $\F{S}$ is Fourier components with respect to $u$ and $S$,
respectively.

Suppose a homogeneous equation relating Eq.(\ref{eq:basic_eq}),
\begin{equation}
  ({\cal D}_t + {\cal K}) \F{u}(\veck,t) = 
  {\cal L}(D) \F{u}(\veck,t) = 0, \label{eq:homo}
\end{equation}
where $D = d/dt$.
The general solution of Eq.(\ref{eq:homo}) is obtained as,
\begin{equation}
  \F{u}(\veck,t) = \sum_j \sum_{l=0}^{m-1} c_{lj}t^l e^{\lambda_j t},
\end{equation}
where $m$ is the multiplicity of the solution $\lambda_j$, $\lambda_j$
is the solution of ${\cal L}(\lambda) = 0$ which is the eigen equation
of Eq.(\ref{eq:homo}).

Let $m$ be 1, which means there is no multiplicity we have simplified
general solution of Eq.(\ref{eq:homo}).
\begin{equation}
  \F{u}(\veck,t) = \sum_j c_j e^{\lambda_j t}.
  \label{eq:gen_sol}
\end{equation}

To obtain explicit scheme for the numerical solution of
Eq.(\ref{eq:basic_eq}) we should have a recurrence formulation with
respect to time step $t_n$.
The explicit scheme is obtained directly from Eq.(\ref{eq:gen_sol}) as,
\begin{eqnarray}
  \F{u}(\veck,t_n) & = & \sum_j c_j e^{\lambda_jt_n} \nonumber \\
  & = & \sum_j c_j e^{\lambda_j(t_{n-1} + \Delta t_n)} \nonumber \\
  & = & \sum_j c_j e^{\lambda_jt_{n-1}} e^{ \lambda_j\Delta t_n}, \nonumber
\end{eqnarray}
where
\begin{equation}
  \Delta t_n = t_n - t_{n-1}.
\end{equation}

Here the authors introduce a state variable $F_n^j$ for numerical formulation,
which is defined as, 
\begin{equation}
  F_n^j = c_j e^{\lambda_jt_n}.
  \label{eq:homo_f}
\end{equation}
Then the authors have explicit scheme for the Eq.(\ref{eq:homo}),
\begin{eqnarray}
  \F{u}(\veck,t_n) & = & \sum_j F_n^j, \nonumber \\
  & = & \sum_j F^{n-1}_{j} \cdot e^{\lambda_j\Delta t_n}.
  \label{eq:tgm_gen_sol}
\end{eqnarray}

General solution of Eq.(\ref{eq:homo}) have been evaluated with a condition of
no-multiplicity of the eigen values.
The derived formulation is recursive with respect to time step.
This means $\F{u}(\veck,t_n)$ can be obtained from the state variable of
one time step before.

Now the authors consider the general solution of Eq.(\ref{eq:homo}).
The authors introduce an impulse function in order to obtain the general
numerical solution of the inhomogeneous equation (\ref{eq:basic_eq}), 
\begin{equation}
	I(t,t_n ;\Delta t_n) = \label{eq:impulse}
	\left\{
	\begin{array}{ll}
		0 & ( t \leq t_{n-1}) \\
		1 & ( t_{n-1} < t < t_{n} )  \\
		0 & ( t_{n} \leq t)
	\end{array}
	\right. .
\end{equation}

\begin{figure}
  \begin{center}
    \includegraphics{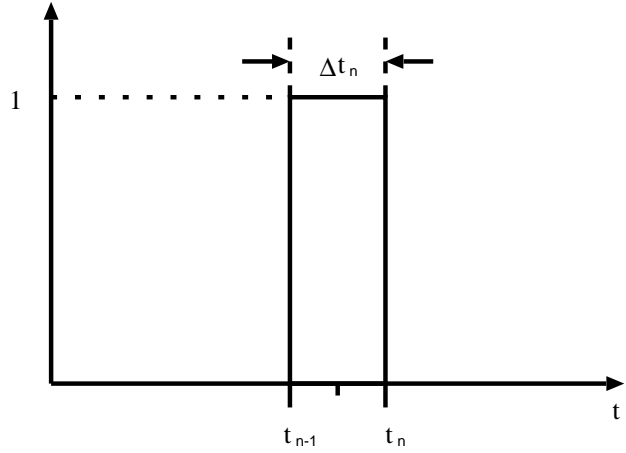}
    \caption{Impulse function for the Green's function. }
    \label{fig:impulse}
  \end{center}
\end{figure}
\begin{figure}[h]
  \begin{center}
    \includegraphics{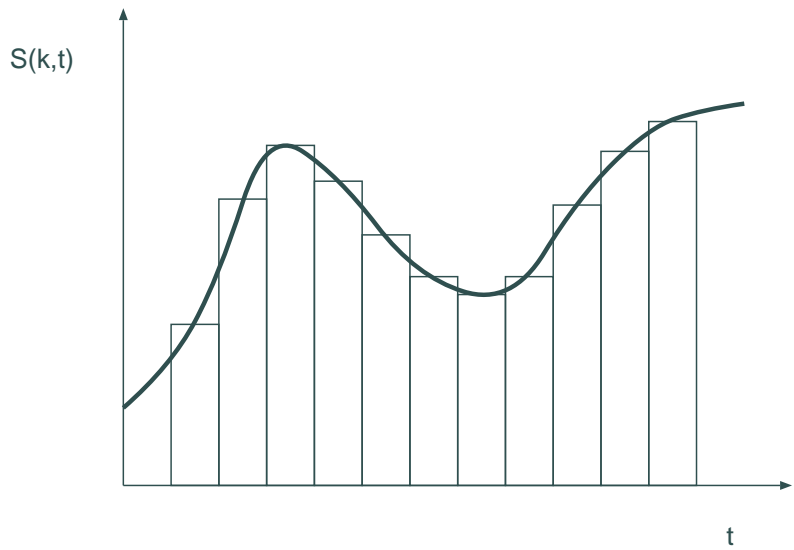}
    \caption{Approximation with the superposition of the Impulse functions.}
    \label{fig:green}
  \end{center}
\end{figure}

Suppose the Green's function $G_n(t)$ of Eq.(\ref{eq:basic_eq}) with
respect to the impulse function which satisfies the equation;
\begin{equation}
  {\cal D}_t G_n(t) + {\cal K} G_n(t) = I(t;t_n,\Delta t).
  \label{eq:f_green}
\end{equation}

A certain source function $S(\veck,t)$ can be approximated with
superposition of the impulse functions as,
\begin{equation}
  \F{S}(\veck,t) \approx \sum_{n=0}^N \F{S}_n(\veck,t_{n-1/2}) I(t;t_n,\Delta t),
  \label{eq:sn}
\end{equation}
where the time series $\F{S}_n(\veck,t_{n-1/2})$ is given.

Equation (\ref{eq:f_green}) is a linear equation so a superposition of
the solutions is also its solution.
From Eqs.(\ref{eq:sn}),(\ref{eq:f_green}) we obtain a superposed
solution as,
\begin{equation}
  \F{u}(\veck,t) \approx \F{u}_N(\veck,t) = \sum_{n=0}^N \F{S}_n(\veck,t_n)
  G_n(t),
  \label{eq:super}
\end{equation}
where $t \geq t_N$.

The authors split summation in Eq.(\ref{eq:super}) into two parts to
obtain an explicit scheme for the Eq.(\ref{eq:basic_eq}),
\begin{equation}
  \F{u}_N(\veck,t) = \sum_{n=0}^{N-1} \F{S}_n(\veck,t_{n-1/2}) G_n(t) 
  + \F{S}_N(\veck,t_{N-1/2}) G_N(t).
  \label{eq:gen_ex}
\end{equation}
If $t \geq t_N$, the first term of the right hand of the
Eq.(\ref{eq:gen_ex}) is the solution of Eq.(\ref{eq:homo}).
So the first term can be described as,
\begin{eqnarray}
  \sum_{n=0}^{N-1} \F{S}_n(\veck,t_{n-1/2}) G_n(t) & = & \sum_j c_j
  e^{\lambda_j t}, \nonumber \\
  & = & \sum_j F^j_{N-1} \cdot e^{\lambda_j (t-t_{N-1})}.
\end{eqnarray}

Now we have the explicit scheme to solve Eq.(\ref{eq:basic_eq}) with the
proposed method in a general form.
\begin{equation}
  \F{u}_N(\veck,t) = \sum_j F^j_{N-1} \cdot e^{\lambda_j (t-t_{n-1})}
  + \F{S}_N(\veck,t_{N-1/2}) G_N(t).
  \label{eq:gen_ex2}
\end{equation}

In the following section we will investigate two concrete examples using the
proposed scheme.

%$\F{u}_n(\veck,t)$ is expected to be described as,
%\begin{equation}
%  \F{u}_n(\veck,t) = \sum_j F^J_{n-1} \cdot e^{\lambda_jt}
%\end{equation}

\section{Formulation Examples}
\label{sec:tgm}
We present two examples of the proposed formulation for solving 
partial differential equations which is derived in the present section.
First one is a wave equation and second is a diffuse equation.
They both have constant coefficients.

\subsection{Wave Equation}
\label{sec:wave}

Here we consider a wave equation of constant coefficients.
\begin{equation}
  \Inv{c^2} \Difd{u(\vecx,t)}{t} - \Laplace u(\vecx,t) = S(\vecx,t),
  \label{eq:wave_equation}
\end{equation}
where $u$ is a physical value, $\vecx$ and $t$ are the coordinates in
space and time respectively, $c$ is the speed of the wave propagation.
$S$ is a source function and is given.
We shows the formulation procedure of Eq.(\ref{eq:wave_equation}) in the
following.

Applying Fourier transform to Eq.(\ref{eq:wave_equation}) gives,
\begin{equation}
  \Inv{c^2} \Difd{\F{u}(\veck,t)}{t} + k^2 \F{u}(\veck,t) =
  \F{S}(\veck,t),
  \label{eq:f_wave_equation}
\end{equation}
where $k^2 = k_x^2 + k_y^2 + k_z^2$.
Now we have an ordinary differential equation instead of partial
differential equation.

Consider the Green's function of Eq.(\ref{eq:f_wave_equation}) against
the impulse function (Eq.(\ref{eq:impulse})), we have an equation,
\begin{equation}
  \Inv{c^2} \Difd{G_n(t)}{t} + k^2 G_n(t) = I(t,t_n;\Delta t_n).
  \label{eq:ft_eq_wave}
\end{equation}
Then the eigen equation can be obtained as,
\begin{equation}
  L(\lambda) = (\lambda^2 + (kc)^2) = 0,
\end{equation}
so we have the eigen solutions.
\begin{equation}
  \lambda = \pm ikc.
\end{equation}

According to Eq.(\ref{eq:gen_ex2}) the solution can be described as,
\begin{eqnarray}
  \F{u}_N(\veck,t) & = & F^+_{N-1} \cdot e^{ikc(t-t_{N-1})} 
  + F^-_{N-1} \cdot e^{-ikc(t-t_{N-1})} \nonumber \\
  & & + \F{S}_N(\veck,t_N)G_N(t),
\end{eqnarray}
where $F_N^+$ and $F_N^-$ are the state variables for $\lambda=ikc$ and
$\lambda=-ikc$, respectively.

Equation (\ref{eq:ft_eq_wave}) can be solved analytically as, 
\begin{eqnarray}
  G_N(t) & = & - i \frac{\sin (kc\Delta t_N/2)}{k^2}
  \left[
  e^{ikc(t-t_N+\Delta t_N/2)} - e^{-ikc(t-t_N+\Delta t_N/2)}
  \right],
  \label{g_func_wave} \\
  & & (t \geq t_N). \nonumber
\end{eqnarray}
Equation (\ref{g_func_wave}) is the Green's function which satisfies 
Eq.(\ref{eq:ft_eq_wave}). 

We redefine the state variable $F_N$ as,
\begin{equation}
  F^{\pm}_N =   F^{\pm}_{N-1} \cdot e^{\pm ikc\Delta t_N} +
  \F{S}(k,t_{N-1/2}) e^{\pm ikc \Delta t_N/2} \sin \frac{kc\Delta t_N}{2}.
\end{equation}

Then we have the formulation by the proposed scheme.
\begin{equation}
  \F{u}(k,t) = - \frac{i}{k^2} \left[ e^{ikc(t-t_N)} F^+_N
    - e^{-ikc(t-t_N)} F^-_N\right] \hspace{5mm} (t \geq t_N).
  \label{eq:tgm_wave}
\end{equation}

\subsection{Diffuse Equation}
\label{sec:diffuse}
Suppose a diffuse equation with constant diffusion coefficient as,
\begin{equation}
  \Dif{u(\vecx,t)}{t} - c \Laplace u(\vecx,t) = S(\vecx,t),
  \label{eq:diffuse}
\end{equation}
where $c$ is the diffusion coefficient, $u(\vecx,t)$ and $S(\vecx,t)$
are functions of space $\vecx$ and time $t$.
Applying Fourier transformation to the Eq.(\ref{eq:diffuse}) to get the
Fourier transformed equation.
\begin{equation}
\Dif{\F{u}(\veck,t)}{t} + ck^2 \F{u}(\veck,t) = \F{S}(\veck,t).
  \label{eq:f_diffuse_eq}
\end{equation}
Same as the procedure deriving the formulation for the wave equation we
have an equation to consider as,
\begin{equation}
  \Dif{G_n(t)}{t} + ck^2 G_n(t) = I(t,t_n ; \Delta t).
  \label{eq:g_eq_dif}
\end{equation}
Equation (\ref{eq:g_eq_dif}) can be solved analytically and then we
have, 
\begin{equation}
  G_N(t) = \Delta t_N \frac{\sinh ck^2\Delta t_N/2}{ck^2\Delta
t_N/2} e^{-ck^2(t-t_N+\Delta t_N/2)}.
  \label{eq:g_sol_dif}
\end{equation}
Same as the formulation of the wave equation we have the explicit form
of the solution as,
\begin{equation}
  \F{u}_N(\veck,t) = F^N(\veck,t_N) e^{-ck^2(t-t_N)},
  \label{eq:tgm_difuse}
\end{equation}
where
\begin{eqnarray}
  F^N(\veck) & = & F^{N-1}(\veck) e^{-ck^2\Delta t_N} \nonumber \\
  & + & \Delta t \frac{\sinh ck^2 \Delta t_N/2}{ck^2 \Delta t_N/2}
  e^{-ck^2\Delta t_N/2} \F{S}(\veck,t_{N-1/2}). \nonumber
\end{eqnarray}

In the next section we will verify the effectiveness of the obtained
formulations by applying to a simple one dimensional problems.
\section{Numerical Experiments}
\label{sec:numerical}
A numerical experiments by the proposed scheme are shown here.
For reader's understanding simple one dimensional examples are shown but
the scheme can be applied to multi-dimension straightforwardly.

The source term is given and the same source function is used for both
examples.
The source function adopted here is a Gaussian distribution in space and
sinusoidal in time.
\begin{equation}
  S(x,t) =  
  \left\{
    \begin{array}{ll}
      0 & ( t < 0 ) \\
      e^{-\frac{(\vecx-\vecx_0)^2}{2\Delta x^2}} \sin \omega_0 t & ( t
      \geq 0 )
    \end{array}
    \right. ,
\end{equation}
where $x_0, \Delta x$ are the centre and the half-value width of the
Gaussian distribution, respectively.

\begin{figure}[bt]
  \begin{center}
    \subfigure[Spatial distribution of the source function.]{\includegraphics[scale=.5]{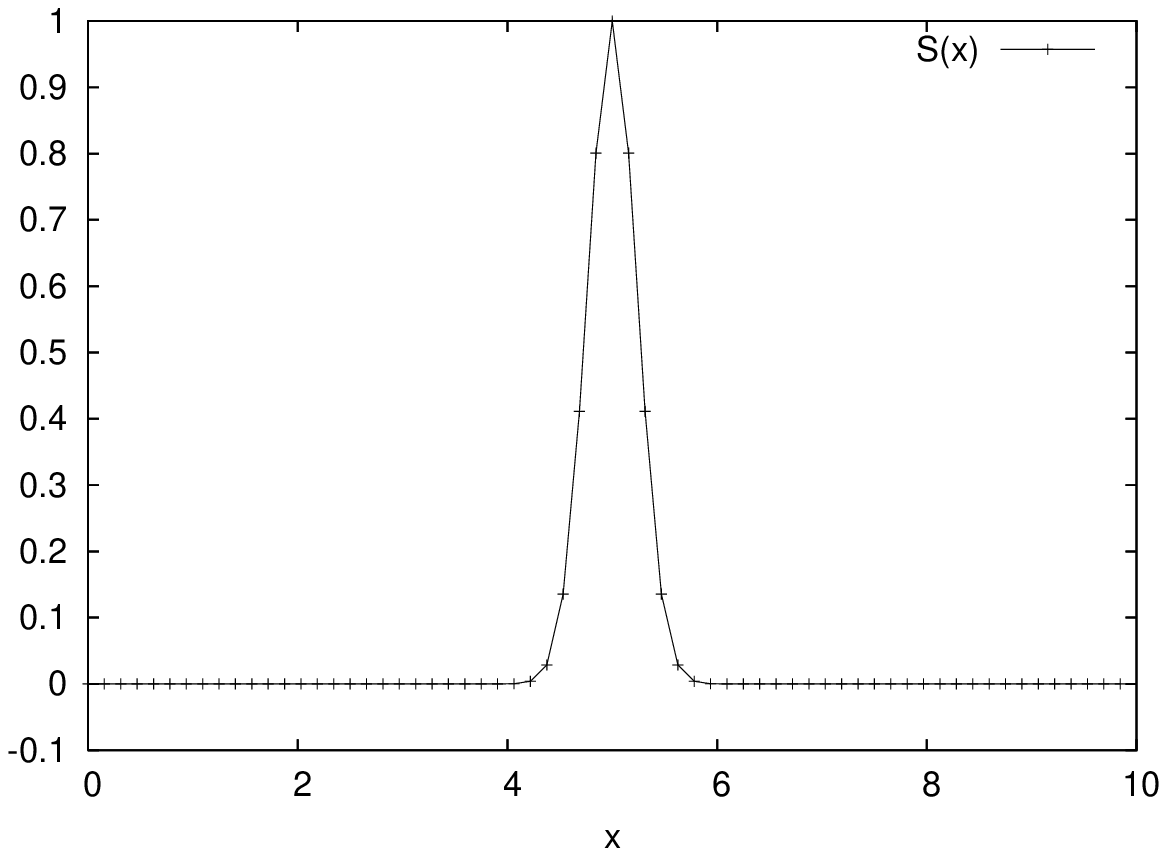} 
      \label{fig:dist}}
    \subfigure[Temporal variation of the source function.]{\includegraphics[scale=.5]{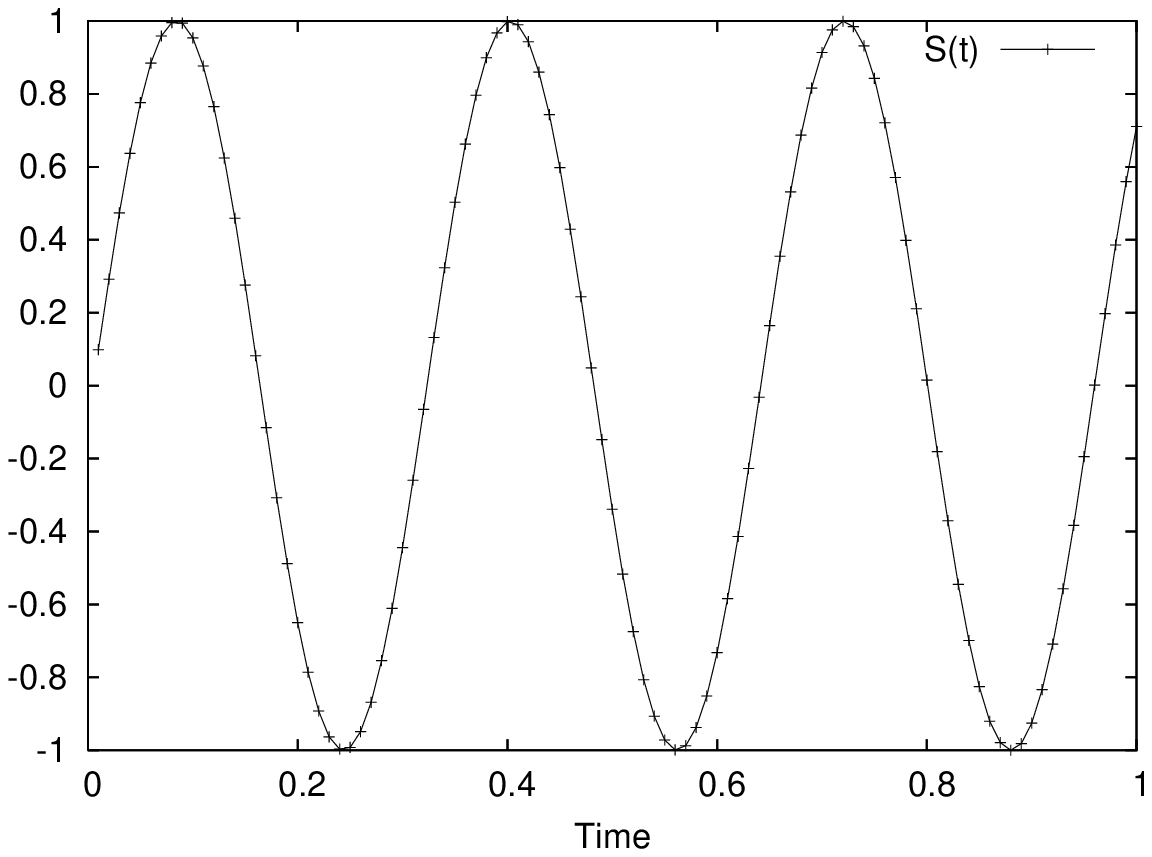} 
      \label{fig:tempo}}
  \end{center}
  \caption{The source function.}
\end{figure}

The Fourier transformed Gaussian function $\F{S}(\veck)$ is analytically
evaluated as,
\begin{equation}
  \F{S}(k) = (2\pi\Delta x^2)^{1/2} e^{-k^2\Delta x^2/2} e^{-i\veck\vecx_0}.
\end{equation}
Figure \ref{fig:dist} and \ref{fig:tempo} shows the distribution of the function
in space and time, respectively.
The common data in both wave and diffuse equations is shown in Table
\ref{tab:common}.
\begin{table}[tb]
  \caption{Common Parameters for the Simulation}
  \label{tab:common}
  \begin{tabular}{c|c|c} \hline
    Grid Number & Length of the system & Grid Width $\Delta x$ 
    \\ \hline \hline
    64 & 10.0 & 10.0/63 \\ \hline
  \end{tabular}
\end{table}

The proposed method is denoted as transient Green method (hereafter TGM)
for convenience.

\subsection{Wave Equation}
\label{sec:num-wave}
A simple one dimensional numerical simulation is performed in order to
investigate the characteristics of the TGM.
The basic equation to solve is Eq. (\ref{eq:f_wave_equation}) and finite
difference solution is also calculated for comparison.

The analytical solution of the basic equation can be obtained and we have,
\begin{equation}
  \F{u}_E (k,t) = \frac{e^{- k^2 \delta_x^2/2 -ikx_0}}
  {\omega_0^2-(ck)^2 } 
  \left[ 
    \frac{\omega}{ck}\sin ckt - \sin\omega_0t \right]. \label{eq:exact_wave}
\end{equation}
The particular parameters for the wave equation is listed in
Table \ref{tab:wave}.
\begin{table}[tb]
  \caption{Parameters for wave equation.}
  \label{tab:wave}
  \begin{tabular}{c|c|c} \hline
    Wave Speed ($c$) & Time Step Width $\Delta t$ & Frequency ($\omega/2\pi$)\\ \hline \hline
    1.0 \hfill &  0.01 \hfill & 3.14159 \\ \hline
  \end{tabular}
\end{table}

\begin{figure}[hbtp]
  \begin{center}
    \includegraphics{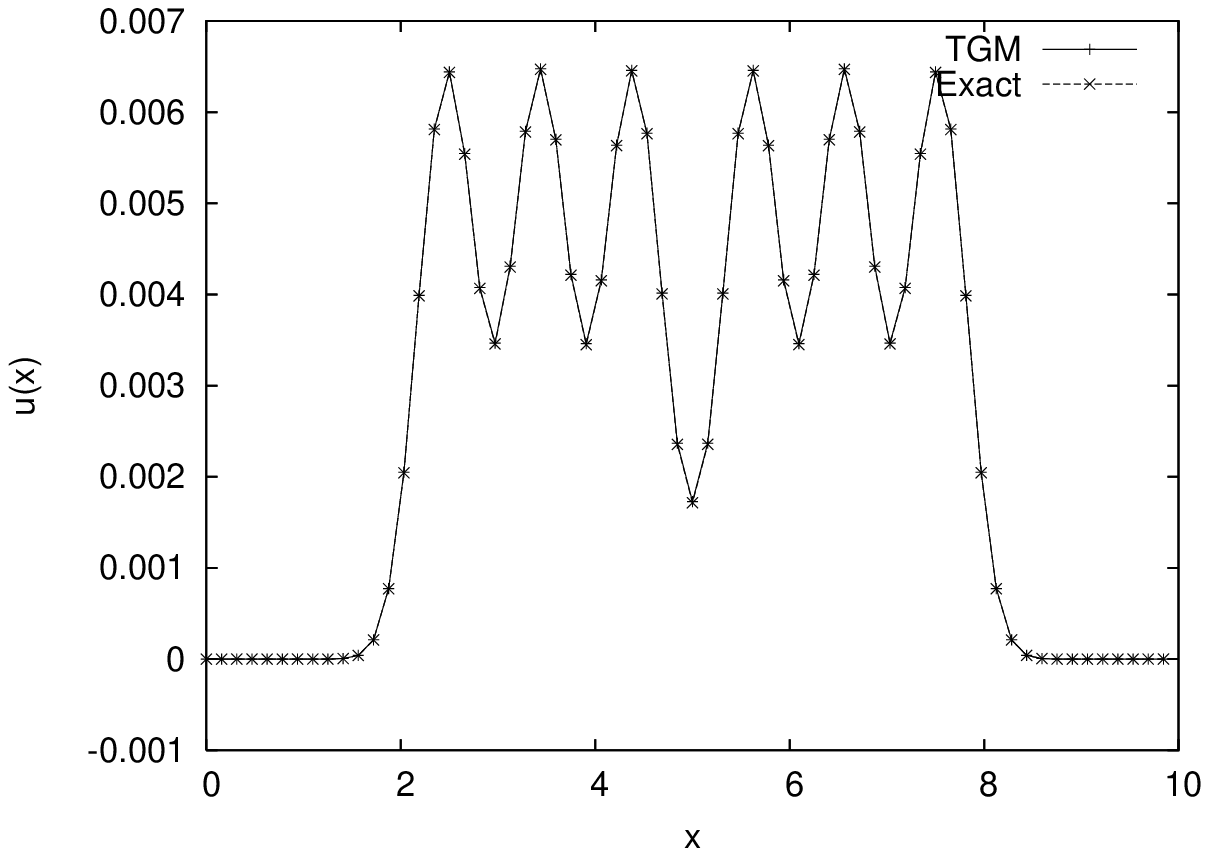}
    \caption{Snapshot of $u_N(\vecx,t)$ and $u_E(\vecx,t)$ at $t=1.0$.} 
    \label{fig:wave_snap}
  \end{center}
\end{figure}
The snapshot at $t=1.0$ is shown in Fig.\ref{fig:wave_snap}.
Although time advancing is performed in $\veck$ space, the figure shows
the real space distribution by inverse Fourier transformation to help
the reader's understanding.
Fast Fourier Transform is used to calculate.
The result shows the good agreement with TGM and exact solution.

\begin{figure}
  \begin{center}
    \includegraphics{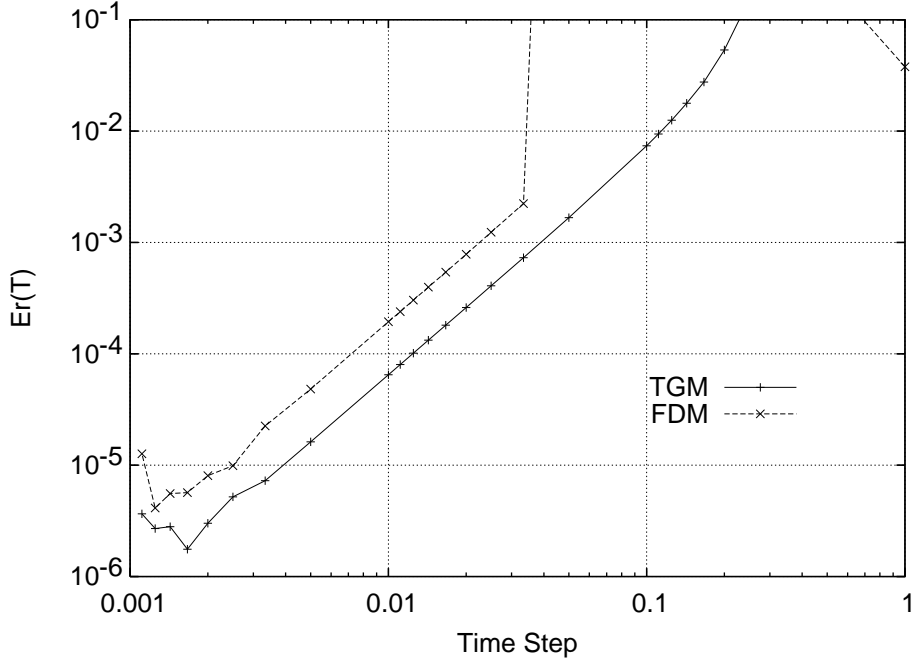}
    \caption{Standard deviation of TGM and FDM from the exact solution
      as a function of $\Delta t$.}
    \label{fig:wave_err}
  \end{center}
\end{figure}
Figure \ref{fig:wave_err} shows the standard deviation from the exact
solutions in $\veck$ space as a function of time step width ($\Delta t$).
\begin{equation}
  Er(t;\Delta t) = \sqrt{ \sum_k ( u(\veck,t;\Delta t) - u_E(\veck,t) )^2
  }.
  \label{eq:error_est}
\end{equation}
In Fig.\ref{fig:wave_err} the FDM and TGM show error with ${\cal
  O}(\Delta t^2)$.
When the time step width $\Delta t$ becomes big enough to violate the
CFL condition the result of FDM shows the numerical instability.
While TGM doesn't show any numerical instability.
It seems to be always stable irrespective of any time step width.
\subsection{Diffuse Equation}
\label{sec:num-heat}
Same as Sec. \ref{sec:num-wave}, numerical simulation is
performed with TGM and FDM for solving a simple one dimensional diffuse
equation problem.
The equation to solve is a diffuse equation
(Eq.(\ref{eq:f_diffuse_eq})). 
The analytical exact solution is obtained for this equation, which is:
\begin{eqnarray}
  \F{u}_E (k,t) & = & \frac{e^{- k^2 \delta_x^2/2 -ikx_0}}{(ck^2)^2 +
\omega_0^2}   \label{model-solution} \\ 
& \times & \left[ \omega_0 e^{-ck^2t} + ck^2
\sin \omega_0t - \omega_0 \cos\omega_0t \right]. \nonumber 
\end{eqnarray}
The particular parameters adopted in the experiment is listed in Table
\ref{tab:diffuse}. 
\begin{table}[hb]
  \caption{Parameters particularly used for Diffuse Equation.}
  \label{tab:diffuse}
  \begin{tabular}{c|c|c} \hline
    Diffuse Coeff. ($c$) & Time Step Width $\Delta t$ & Frequency ($\omega/2\pi$) \\ \hline \hline
    3.0 &  0.001 & 20.0 \\ \hline
  \end{tabular}
\end{table}

\begin{figure}[hbtp]
  \begin{center}
    \includegraphics{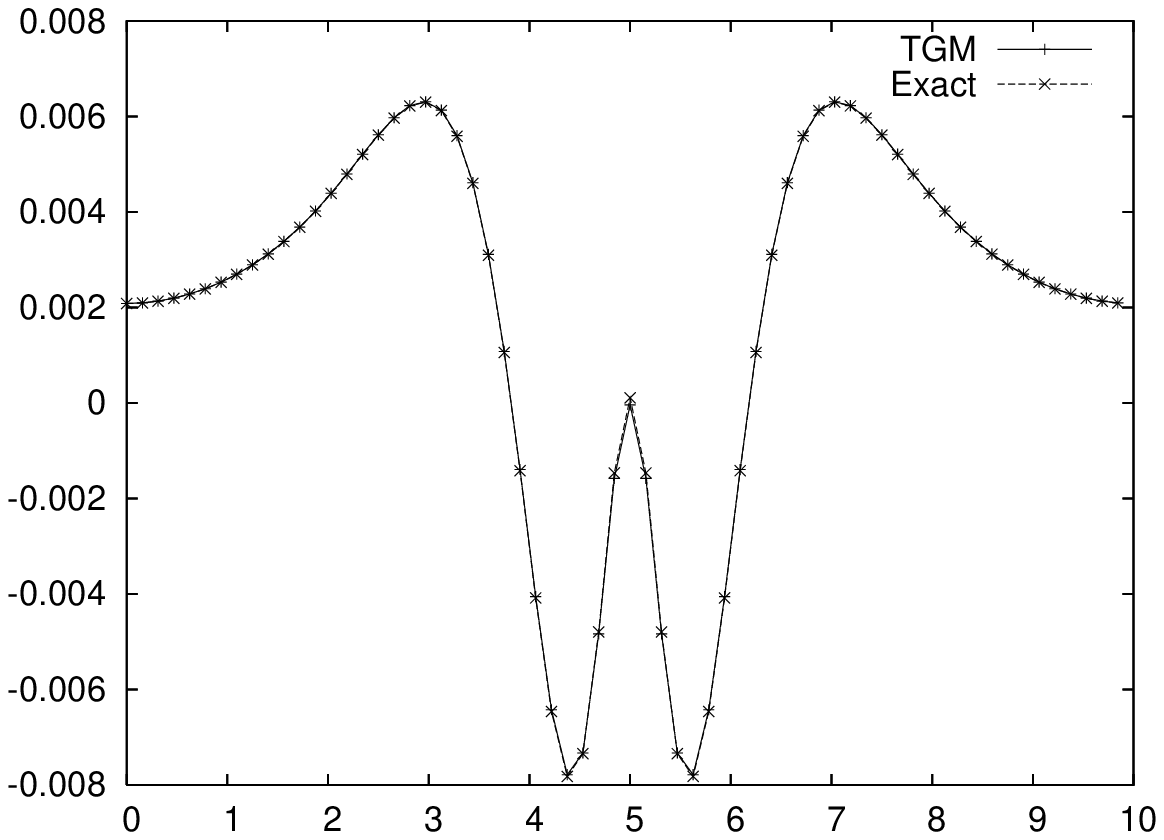}
    \caption{Snapshot of $u_N(\vecx,t)$ and $u_E(\vecx,t)$ at $t=0.1$.} 
    \label{fig:difuse}
  \end{center}
\end{figure}

The snapshot at $t=0.1$ is shown in Fig.\ref{fig:difuse}.
The solution computed by TGM and exact solution indicate good agreement.

\begin{figure}[hbtp]
  \begin{center}
    \includegraphics{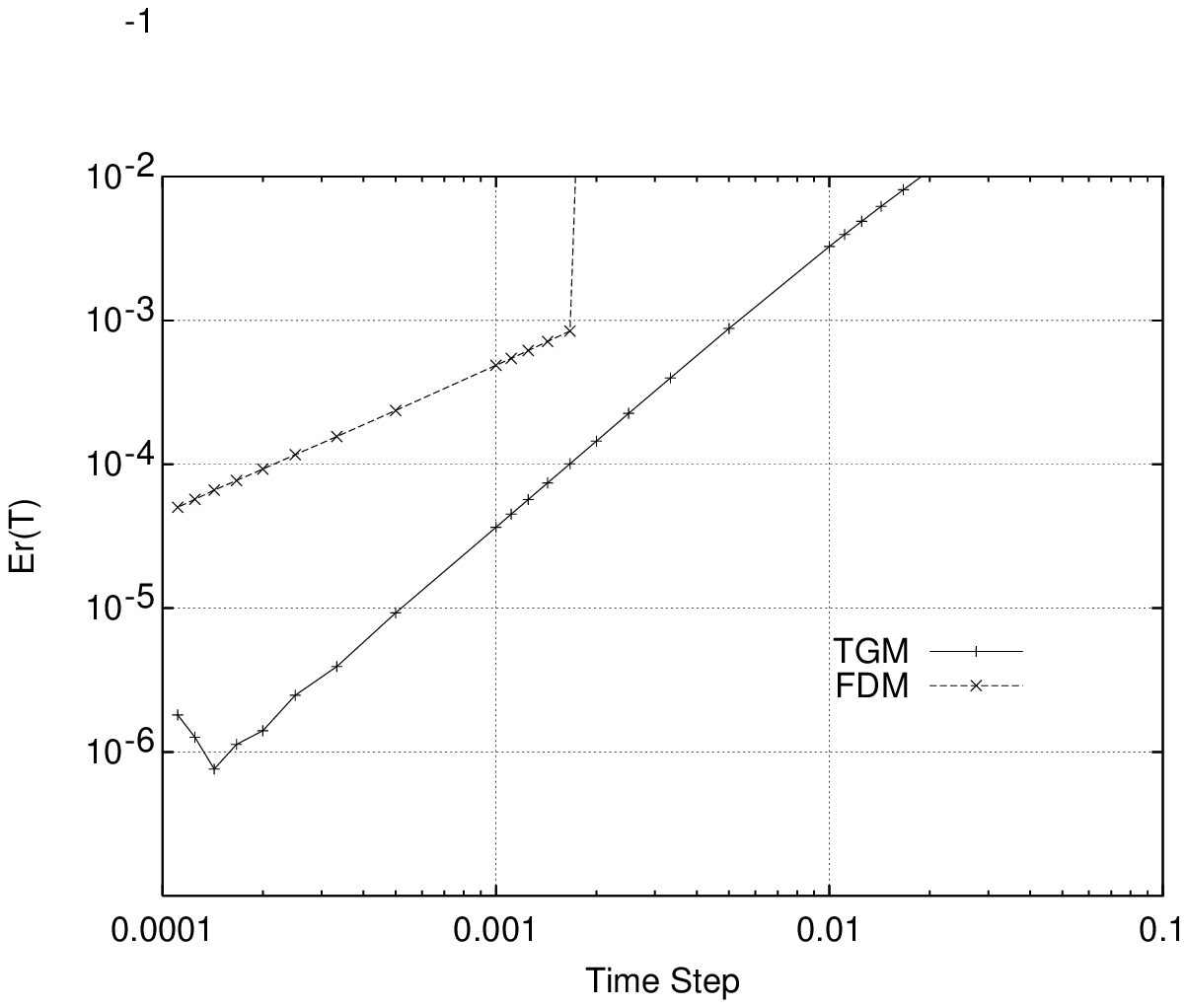}
    \caption{Standard deviation of TGM and FDM from the exact solution
      as a function of $\Delta t$.}
    \label{fig:difuse_err}
  \end{center}
\end{figure}

The errors of the TGM and the FDM are shown in Fig.\ref{fig:difuse_err}.
The error is evaluated form Eq.(\ref{eq:error_est}) same as the case of
wave equation.
Figure \ref{fig:difuse_err} shows numerical instability is occurred with
the violation of CFL condition whereas the results of the TGM doesn't
show any numerical instability irrespective of any time step width.

The FDM is the Euler's scheme so the error is ${\cal O}(\Delta t)$.
According to the Fig.\ref{fig:difuse_err}, TGM has accuracy of ${\cal
  O}(\Delta t^2)$. 

\section{Discussion}
\label{sec:discuss}

\subsection{Numerical Accuracy}
\label{sec:accuracy}
The numerical error of the proposed scheme potentially contains the
discretization errors in time and space.

The spatial accuracy of TGM is determined by the Fourier transformation.
If the Fast Fourier Transformation is adopted the data can be
transformed exactly.
This means the transformed data contains no truncation error essentially.

The temporal accuracy of the proposed method depend upon the accuracy of
the integration method of the source term $S(\veck,t)$.
In the presented scheme here
(Eqs.(\ref{eq:tgm_wave}),(\ref{eq:tgm_difuse})) rectangle mid-point rule
is adopted for time integration (Eq.(\ref{eq:sn})).
So discretization error of Eqs.(\ref{eq:tgm_wave}),(\ref{eq:tgm_difuse})
is ${\cal O}(\Delta t^2)$.
If you adopt more accurate integration scheme the accuracy of the
proposed scheme also becomes more accurate.

\subsection{Numerical Stability}
\label{sec:neumann}
The numerical stability of the proposed method is discussed in the
section. 

As for the finite difference method the explicit scheme is unstable
unless the time step width $\Delta t$ is small enough to satisfy a
condition. 
The criteria is well known as Courant Friedrichs and Lewy (CFL) condition.
For instance CFL condition of the simple explicit scheme for diffuse
equation (\ref{eq:diffuse}) is; 
\begin{equation}
  \Delta t < \frac{2\Delta x^2}{c}. \nonumber
\end{equation}
If the condition is violated the simulation result may be far from the
exact solution which you want to obtain.
The unstableness comes from the difficulty in approximating the solution
with first or second order polynomials.
While the proposed method is derived from the superposition of the
Green's function which are analytical solutions.

Now we consider the numerical stability of the TGM solution for
Eq.(\ref{eq:homo}).
The TGM formulation is given by Eq.(\ref{eq:tgm_gen_sol}).
The general solution of Eq.(\ref{eq:homo}) is already given by
Eq.(\ref{eq:gen_sol}).
The state variable $F_n^j$ also given by Eq.(\ref{eq:homo_f}).

When $F^j_{n-1}$ satisfy the equation:
\begin{equation}
  F^j_{n-1} = c_j e^{\lambda_j t_{n-1}}. \nonumber
\end{equation}
The $F^j_{n}$ can be in TGM formulation as,
\begin{eqnarray*}
  F^j_n & = & F^j_{n-1} \cdot e^{\lambda_j \Delta t_n}, \\
  & = & c_j e^{\lambda_j t_{n-1}} \cdot e^{\lambda_j \Delta t_n}, \\
  & = & c_j e^{\lambda_j t_{n}}.
\end{eqnarray*}
So we find the TGM solution is exactly identical with the analytical
solution.
This is always true irrespective of any time step width $\Delta t$.
This result indicates the TGM formulation is always stable with any time
step width as is shown in section \ref{sec:numerical}

%The stability with the time step for wave and diffuse equations are
%discussed. 

\section{Conclusion}
A numerical method to compute explicitly for transient solution of linear
partial differential inhomogeneous equation with constant
coefficients is proposed.
It is shown that the proposed method gives the transient solution from
the state variables at the one time step before explicitly and always
numerically stable irrespective of any time step width.


\begin{thebibliography}{10}
\bibitem{cfl} R. Courant, K.O. Friedrichs, and H. Lewy, Math.
Ann., 100, 1928, 32.
\bibitem{anderson} D.A. Anderson, J.C.Tannehill, and R.H.Pletcher,
  Computational Fluid Mechanics and Heat Transfer, McGraw-Hill, 1984,
  p108. 
\bibitem{habe:ipsj} Hiroshi Abe, Transaction on Information
    Processing Society of Japan, 1992, p.1006.
\bibitem{habe:cefc93} Hiroshi Abe, Digests of the Fifth
  Biennial IEEE Conference on Electromagnetic Field Computation,
  Harvey Mudd College, Claremont, CA, USA, 1992, MP39.
\bibitem{habe:jsiam93} Hiroshi Abe, Annual Meeting of Japan Society
    of Industry and Applied Mathematics, 1993.
\end{thebibliography}
\end{document}